\documentclass{ws-jktr}

\usepackage{amscd,graphics}

\newcommand{\Z}{{\mathbb Z}}

\newcommand{\ph}{{\varphi}}
\newcommand{\e}{{\varepsilon}}

\newcommand{\al}{{\alpha}}

\newcommand{\de}{{\delta}}
\newcommand{\la}{{\lambda}}
\newcommand{\ga}{{\gamma}}

\newcommand{\te}{{\theta}}

\newcommand{\Om}{{\Omega}}
\newcommand{\io}{{\iota}}
\newcommand{\ti}{\tilde}
\newcommand{\bu}{\bullet}

\newcommand{\pr}{{\mathrm{pr}}}
\newcommand{\inp}{{\mathrm{Int}\,}}

\newcommand{\es}{{\varnothing}}
\newcommand{\pt}{{ \{\mathrm{pt.}\} }}
\newcommand{\lk}{{ \mathrm{lk} }}

\newcommand{\bd}{{\partial}}

\newcommand{\pp}[1]{\langle #1 \rangle}
\newcommand{\jl}[1]{ \{ #1 \} }
\newcommand{\ov}[1]{ \overline{#1} }
\newcommand{\vmu}{ \bm{\mu} }
\newcommand{\vla}{ \bm{\la} }
\newcommand{\vta}{ \bm{\tau} }
\newcommand{\va}{ {\mathbf a} }
\newcommand{\ved}{ {\mathbf e} }

\catchline{}{}{}{}{}

\title{
PERIPHERALLY SPECIFIED HOMOMORPHS OF LINK GROUPS }

\author{V. KURLIN}

\address{
 Department of Mathematical Sciences,\\
 University of Liverpool,\\
 Liverpool L69 7ZL, United Kingdom\\
 kurlin@liv.ac.uk}

\author{D. LINES}

\address{
 Institut de Math\'ematiques de Bourgogne,\\
 Universit\'e de Bourgogne,\\
 BP 47870, 21078 Dijon cedex, France\\
 dlines@u-bourgogne.fr }

\begin{document}

\maketitle

\begin{abstract}
Johnson and Livingston have characterized peripheral structures
 in homomorphs of knot groups.
We extend their approach to the case of links.
The main result is an algebraic characterization of
 all possible peripheral structures in certain
 homomorphic images of link groups.
\end{abstract}

\keywords{Link, link group, longitude, meridian,
 Pontryagin product, Johnson-Livingston product}

\ccode{Mathematics Subject Classification 2000: 57M25, 57M05}


\section{Introduction}


\subsection{Motivation and summary of results}

Groups that can appear as the image under
 a surjective homomorphism of the group of a knot
 have been investigated by various authors,
 see for instance \cite{3,4,5}.
Given such a homomorphic image $G$, it is of interest
 to characterize the subgroup which is the image of
 the peripheral subgroup of the knot.
D.~Johnson and C.~Livingston \cite{5} have given
 necessary and sufficient conditions on elements
 $\mu$ and $\la$ of $G$ for them to be the image
 of the meridian, respectively the preferred longitude
 of the knot.
These conditions involve a Pontryagin product
 $\pp{\mu,\la}$ in the homology group $H_2(G)$
 and a Johnson-Livingston product $\jl{\mu,\la}$
 in a quotient of $H_3(G/G')$.
\smallskip

We extend their method to the case of $r$-component links.
In this context we consider systems of elements
 $\vmu=(\mu_1,\dots,\mu_r)$ and $\vla=(\la_1,\dots,\la_r)$
 which are the images of the meridians, respectively
 the longitudes of the components of the link.
We show that provided the $\mu_i$ are conjugate in $G$
 one can define an extended Johnson-Livingston product
 $\jl{\vmu,\vla}$ and give necessary and sufficient conditions
 for the realizability of these systems.


\subsection{Preliminary definitions}

\begin{definition}
 (\emph{knots}, \emph{ribbon links}, \emph{ambient isotopy})
\smallskip

\noindent
{(a)}
\emph{A link} $L$ is a topologically flat embedding of
 disjoint oriented circles in $S^3$.
The $i$th circle is called \emph{the $\mathrm{i}$th component} of
 the link $L$ and is denoted by $L_i$.
If $r=1$, then the link $L$ is \emph{a knot}.
A link $L\subset S^3$ is \emph{ribbon}, if $L$ bounds a
 disjoint union of disks $\sqcup_{i=1}^r D_i^2$ immersed into
 $S^3$, whose singularities are always as on Fig.~1.
\smallskip

\noindent
{(b)}
Two links $L$ and $L'$ are \emph{equivalent}, if
 there is an orientation preserving self homeomorphism of $S^3$
 sending $L_i$ to $L'_i$ for $i=1,\dots,r$ and respecting
 the orientations of the components.
Links will be studied up to equivalence.
\end{definition}

\begin{definition}
 (\emph{meridians}, \emph{longitudes},
 \emph{preferred systems} of longitudes)
\smallskip

\noindent
{(a)}
Let $L=\cup_{i=1}^r L_i\subset S^3$ be an $r$-component link,
 $T(L_i)$ be a sufficiently small tubular neighbourhood of $L_i$,
 and $p_i$ a point on boundary $\bd T(L_i)$, $i=1,\dots,r$.
\emph{A meridian} $m_i$ of the component $L_i$ is an oriented
 simple closed curve $m_i\subset\bd T(L_i)$ that bounds
 inside $T(L_i)$ a 2-dimensional disk intersecting $L_i$
 in a single point with positive sign.
The homotopy class of $m_i$ is unique in $\pi_1(\bd T(L_i),p_i)$
 and called \emph{the meridian} of $L_i$.
\smallskip

\noindent
{(b)}
\emph{A longitude} of the component $L_i$ is an oriented curve
 $l_i\subset\bd T(L_i)$ passing through $p_i$ isotopic to
 $L_i$ inside $T(L_i)$.
A longitude of $L_i$ is \emph{preferred} and is denoted by $\bar l_i$,
 if there is an oriented surface $F_i\subset S^3-\inp T(L_i)$ with
 $\bd F_i=\bar l_i$.
The homotopy class of $\bar l_i$ is unique in $\pi_1(\bd T(L_i),p_i)$
 and called \emph{the preferred longitude}.
\smallskip

\noindent
{(c)}
Denote by $T(L)$ the disjoint union $\sqcup_{i=1}^r T(L_i)$
 of sufficiently small tubular neighbourhoods of $L_1,\dots,L_r$.
A system of curves $\ov{(l_1,\ldots,l_r)}\subset S^3-L$
 is \emph{a preferred system of longitudes for $L$},
 if the curve $l_i$ is a longitude of the component $L_i$
 for each $i=1,\ldots,r$,
 and the union $l_1\cup\ldots\cup l_r$ is the boundary
 of an oriented surface $F\subset S^3-\inp T(L)$.
\end{definition}

Defnition~1.2c provides a natural extension
 to links of the notion of a preferred longitude for knots.
Note that $\ov{(l_1,\ldots,l_r)}\neq(\bar l_1,\dots,\bar l_r)$ in general.
The definition of preferred longitudes and systems
 will be reformulated in Lemma~2.2.
\smallskip

We fix a system of arcs $\ga_i$, $i=2,\dots,r$, properly embedded
 in $S^3-\inp T(L)$ joining $p_1$ to $p_i$ and intersecting only
 at $p_1$.
We define $\ga_1$ to be the constant path at $p_1$.
To each oriented simple closed curve $c$ in $\bd T(L_i)$ passing
 through $p_i$, $i=1,\dots,r$ we can associate the homotopy class
 of the loop $\ga_i\circ c\circ\ga_i^{-1}$.
We still denote this homotopy class by $c$ and consider it as
 an element of $\pi_1(S^3-L):=\pi_1(S^3-L,p_1)$.
This holds in particular for the meridians and longitudes of $L$.
\smallskip

Further we are dealing with an abstract finitely generated group $G$
 and surjective homomorphisms $\pi_1(S^3-L)\to G$,
 where $L$ is a link with $r$ components.

\begin{definition}
 (\emph{a meridional system} for $G$ and
 \emph{realizable systems} for $(G,\vmu)$)
\smallskip

\noindent
{(a)}
Let $G$ be a finitely generated group.
Denote by $G^r$ the direct product $G\times\cdots\times G$
 ($r$ times).
A system of elements $\bm{\mu}=(\mu_1,\ldots,\mu_r)\in G^r$
 is called \emph{meridional} for $G$, if the group $G$
 is generated by finitely many conjugates of
 the elements $\mu_1,\ldots,\mu_r$.
\smallskip

\noindent
{(b)}
Suppose that a group $G$ has a meridional system
 $\vmu=(\mu_1,\ldots,\mu_r)\in G^r$.
A system of elements $\vla=(\la_1,\ldots,\la_r)\in G^r$ is called
 \emph{weakly realizable} for the pair $(G,\vmu)$, if
 there exists a link $L=\cup_{i=1}^r L_i\subset S^3$ with
 a surjective homomorphism $\rho:\pi_1(S^3-L)\to G$ such that
 $\rho(m_i)=\mu_i$, $\rho(l_i)=\la_i$ for each $i=1,\ldots,r$, where
 $m_i$ is the meridian of the component $L_i$,
 and $l_i$ is a longitude of $L_i$.
\smallskip

\noindent
{(c)}
A weakly realizable system $\vla=(\la_1,\ldots,\la_r)\in G^r$ is
 \emph{realizable} for the pair $(G,\vmu)$, if
 $(l_1,\ldots,l_r)$ from {(b)} is
 a preferred system of longitudes for the link $L$.
Denote by $R(G,\vmu)\subset G^r$ the set of all realizable
 systems $\vla\in G^r$ for $(G,\vmu)$.
\smallskip

\noindent
{(d)}
Let $G$ be a group with a meridional system $\vmu\in G^r$.
Denote by $G'$ the commutator subgroup of $G$.
Let $\pr:G\to G/G'$ be the quotient map.
Denote by $[\mu_1],\dots,[\mu_r]$ the images of the meridians
 $\mu_1,\dots,\mu_r\in G$ in \emph{the abelianization} $G/G'$.
\end{definition}


\subsection{Main results}

\noindent
Let $G$ be a finitely generated group.
Suppose that there exist an $r$-component link $L\subset S^3$
 and a surjective homomorphism $\rho:\pi_1(S^3-L)\to G$.
For any link $L\subset S^3$, the group $\pi_1(S^3-L)$ has a well-known
 Wirtinger presentation \cite{2} and hence a meridional system.
Then $G$ also has a meridional system obtained by selecting
 one Wirtinger generator $m_i$ for each component $L_i$.
Their images $\mu_i=\rho(m_i)$ form a meridional system
 $\vmu=(\mu_1,\ldots,\mu_r)\in G^r$.
\smallskip

The converse was proved by Gonzalez-Acuna for links \cite{3}
 by using a 4-dimensional technique.
We shall need the fact that every meridional system is
 realized by a link where all the linking numbers are zero,
 so we extend to ribbon links the simple geometric proof
 given for knots in \cite{4}, Proposition~2.3.
This result shows that the set $R(G,\vmu)$ is not empty.

\begin{theorem}
Let a group $G$ have a meridional system
 $\vmu=(\mu_1,\ldots,\mu_r)\in G^r$.
\smallskip

\noindent
A system $\vla=(\la_1,\ldots,\la_r)\in G^r$
 is \emph{weakly realizable} for $(G,\vmu)$ if and only if
\smallskip

\noindent
(i) the element $\la_i$ commutes with $\mu_i$ for each
 $i=1,\ldots,r$;
\smallskip

\noindent
(ii) the sum of the Pontryagin products
 $\sum\limits_{i=1}^r \pp{\mu_i,\la_i}$ vanishes in the group $H_2(G)$.
\end{theorem}

Let $Q(G)$ denote the quotient group $H_3(G/G')/pr_*(H_3(G))$.

\begin{theorem}
Let a group $G$ have a meridional system
 $\vmu=(\mu_1,\ldots,\mu_r)\in G^r$ such that
 each $\mu_i$ is conjugate to $\mu_1$, $i=2,\dots,r$.
A system of elements $\vla=(\la_1,\ldots,\la_r)\in G^r$
 is \emph{realizable} for the pair $(G,\vmu)$ if and only if
 Conditions~(i),(ii) of Theorem~1.4 hold and
\smallskip

\noindent
(iii) $\vla\in(G')^r$;
\smallskip

\noindent
(iv) the extended Johnson-Livingston product
 $\{\vmu,\vla\}=0$ vanishes in $Q(G)$.
\end{theorem}

To extend the Johnson-Livingston method we use multi-connected sums
 of various geometric objects such as links, surfaces and manifolds,
 see Definitions~4.5 and 4.10.
The key points of the proof are the well-definedness and
 additivity of the extended Johnson-Livingston product,
 see Theorem~4.17.


\subsection{Organization of the paper}

\noindent
The realizability of meridional systems is proved in section~2,
 Proposition~2.3.
Section~3 contains Definition~3.2 of the Pontryagin product and
 the proof of Theorem~1.4.
In section~4 we introduce Johnson-Livingston products in
 Definitions~4.1, 4.4 and prove their well-definedness.
The proof of Theorem~1.5 will be finished in section~5.
We give examples of applications of Theorem~1.5 in section~6.


\section{Preferred Longitudes and Meridional Systems}

Subsection~2.1 discusses preferred systems of longitudes.
In subsection~2.2 the realizability of meridional systems
 is proved using geometric operations on algebraic representations.


\subsection{Preferred systems of longitudes}

\begin{definition}
 (\emph{the linking number} $\lk$,
 \emph{algebraically split links})

\noindent
{(a)}
Let $J$ and $K$ be two disjoint oriented simple closed curves
 in $S^3$.
We denote by $\lk(J,K)$ their linking number.
\smallskip

\noindent
{(b)}
Let $L=\cup_{i=1}^r L_i\subset S^3$ be an oriented $r$-component link.
If all the linking numbers $\lk(L_i,L_j)=0$ for $i,j=1,\dots,r$, $i\neq j$,
 then the link $L$ is called \emph{algebraically split}.
\end{definition}

Recall that preferred longitudes and systems were introduced in
 Definition~1.2.

\begin{lemma}
Let $L=\cup_{i=1}^r L_i\subset S^3$ be an oriented $r$-component link.
Let $[m_i],[l_i]\in H_1(S^3-L)$ be the classes of
 the meridian and a longitude of $L_i$, $i=1,\dots,r$.
The system of curves $(l_1,\ldots,l_r)$
 is preferred for the link $L$ if and only if
 $[l_i]=\sum\limits_{i=1}^{r}\al_{ij}[m_j]$,
 where $\al_{ij}=\lk(L_i,L_j)$, $i\neq j$,
 $\al_{ii}=-\sum_{j\neq i}\lk(L_i,L_j)$, $i,j=1\dots,r$.
\end{lemma}
\begin{proof}
The homology group $H_1(S^3-L)$ is isomorphic to
 $\Z[m_1]\oplus\cdots\oplus\Z[m_r]$, where $m_i$ is the meridian of the
 component $L_i$, $i=1,\dots,r$.
By Definition~1.2c a system $(l_1,\dots,l_r)$ of
 longitudes is preferred, if there is an oriented surface
 $F\subset S^3-\inp T(L)$ with the boundary $\bd F=\sqcup_{i=1}^r l_i$.
\smallskip

Hence in the group $H_1(S^3-L)$ one has $\sum_{i=1}^r[l_i]=0$.
The preferred longitude $\bar l_i$ is the boundary of
 an oriented surface $F_i\subset S^3-\inp T(L_i)$, $i=1,\dots,r$.
The surface $F_i$ can be used to compute
 the linking number $\al_{ij}=\lk(L_i,L_j)$, $i\neq j$.
Then $[\bar l_i]=\sum_{j\neq i}\al_{ij}[m_j]$ in $H_1(S^3-L)$.
\smallskip

By Definition~1.2b any longitude $l_i\subset\bd T(L_i)$ is
 isotopic inside $T(L_i)$ to $L_i$, hence $l_i=m_i^{\al_{ii}}\bar l_i$
 in $\pi_1(\bd T(L_i))$ for some $\al_{ii}\in\Z$, $i=1,\dots,r$.
Equivalently, in the homology group $H_1(S^3-L)$ one gets
 $[l_i]=\al_{ii}[m_i]+[\bar l_i]=\sum_{j=1}^r\al_{ij}[m_j]$.
The condition $\sum_{i=1}^r[l_i]=0$ is equivalent to
 $\al_{ii}=-\sum_{j\neq i}\al_{ij}$ as required.
\smallskip

Conversely, construct a smooth map $f:\sqcup \bd T(L_i)\to S^1$
 such that the restriction to the meridians $m_i$ of $L_i$
 is a degree one map, and $f$ is constant on the curves $l_i$.
The only obstructions to the extension of the map $f$
 to a smooth map $\ti f:S^3-\sqcup_{i=1}^r\inp T(L_i)\to S^1$
 are the conditions $\al_{ii}+\sum_{j\neq i}\al_{ij}=0$,
 $i=1,\dots,r$.
\smallskip

An inverse image of a regular value of the extended map $\ti f$
 gives an oriented surface $F\subset S^3-\sqcup_{i=1}^r\inp T(L_i)$
 such that the boundary $\bd F=l_1\cup\dots\cup l_r$.
Hence the system of curves $(l_1,\dots,l_r)$ is preferred
 for the link $L$.
\end{proof}


\subsection{Realizability of meridional systems}

\begin{proposition}
Suppose that a group $G$ has a meridional system
 $\vmu=(\mu_1,\ldots,\mu_r)\in G^r$.
Then there exists a ribbon $r$-component link $L\subset S^3$
 with a surjective homomorphism $\rho:\pi_1(S^3-L)\to G$
 such that $\rho(m_i)=\mu_i$ for each $i=1,\ldots,r$, where
 $m_i$ is the meridian of the component $L_i$.
\end{proposition}

The proof of Proposition~2.3 is a straightforward generalization
 of Johnson's proof for knots.
We refer the reader to \cite{4} for details.
Recall that to describe a homomorphism from the link group
 to the group $G$, it suffices to label the arcs of a planar
 diagram of the link with elements of $G$ in such a way that,
 at each crossing, the corresponding relation for the labelled
 elements holds in $G$.

\emph{A represented band} is a pair of parallel, oppositely directed
 arcs of the diagram of the link, with no other arcs of the link
 passing between the two arcs of the band and such that
 the two arcs are labelled with the same element of $G$, see
 Fig.~1.

\begin{figure}[h]
\includegraphics{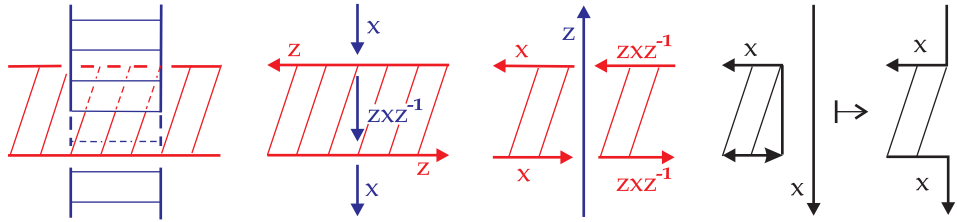}
\caption{
A part of a ribbon link, represented bands and a band connection.}
\end{figure}

\noindent
{\it Proof of Proposition 2.3.}
Let $\mu_1,\dots,\mu_r$ be a meridional system for $G$.
There are positive integers $k_1,\dots,k_r$ and
 $\mu_{ij}\in G$, $i=1,\dots,r$, $j=1,\dots,k_i$
 such that:
\smallskip

$\bu$
 $\mu_{i1}=\mu_i$ for all $i=1,\dots,r$;
\smallskip

$\bu$
 the $\mu_{ij}$ generate the group $G$;
\smallskip

$\bu$
 $\mu_{ij}=W_{ij}\mu_i W_{ij}^{-1}$ for words
 $W_{ij}$ in $\mu_{ij}$, $i=1,\dots,r$, $j=1,\dots,k_i$.
\medskip

Consider the trivial link of $k_1+\cdots+k_r$ components
 and label the components with the $\mu_{ij}$ (see Fig.~2).
This gives a surjective homomorphism from the group of the
 trivial link to $G$.
For each pair of $i,j$ weave a represented band issuing
 from the component labelled with $\mu_{i1}$ according to the
 instructions given by $W_{ij}$.
\smallskip

Perform a band connection (see Fig.~1) between
 the represented band and the component labelled with $\mu_{ij}$.
One obtains in this way a ribbon link $L$ of $r$ components
 and a surjective homomorphism $\rho:\pi_1(S^3-L)\to G$ (see Fig.~2).
\hfill $\square$

\begin{figure}[h]
\includegraphics{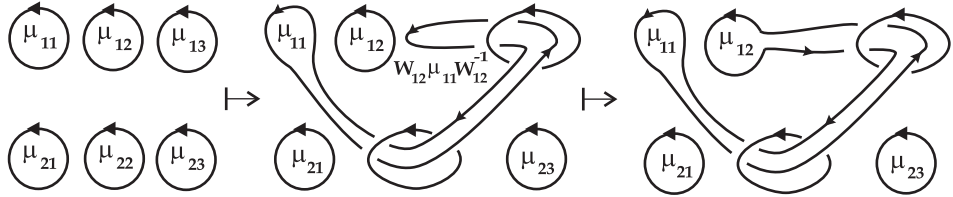}
\caption{
How to construct a ribbon link $L$ with
 a homomorphism $\rho:\pi_1(S^3-L)\to G$?}
\end{figure}


\section{Weakly realizable systems: proof of Theorem~1.4}

In subsection~3.1 we introduce the Pontryagin product.
Subsection 3.2 is devoted to necessity in Theorem~1.4.
Subsection~3.3 contains the proof of sufficiency in Theorem~1.4.


\subsection{Pontryagin product}

\begin{definition}
 (\emph{the $G$-bordism group} $\Om_n(G)$ of
 $n$-dimensional manifolds)
\smallskip

\noindent
{(a)}
Let $M,N$ be two oriented closed
 $n$-dimensional possibly disconnected manifolds and
 let $G$ be a group.
Fix homotopic classes of continuous maps
 $f_M:M\to K(G,1)$, $f_N:N\to K(G,1)$.
The pairs $[M,f_M]$ and $[N,f_N]$ are called
 \emph{$G$-bordant}, i.e. $[M,f_M]=[N,f_N]$,
 if there is an oriented compact connected
 $(n+1)$-dimensional manifold $W$
 with a continuous map $f_W:W\to K(G,1)$
 such that $\bd W=M\cup(-N)$,
 $f_W|_{M}=f_M$, and $f_W|_{N}=f_N$.
Here $(-N)$ is $N$ with the reversed orientation.
\smallskip

\noindent
{(b)}
Classes of $G$-bordant pairs $[M,f_M]$ form
 \emph{the $G$-bordism group $\Om_n(G)$}.
The operation is the disjoint union, the unit element
 is the empty set $\es$ with the empty map $\es\to K(G,1)$.
\smallskip

\noindent
{(c)}
Denote by $H_n(G)$ the $n$th homology of a group $G$
 with integer coefficients \cite{1}.
It is a well-known fact, which can be proved using the
 Atiyah-Hirzebruch spectral sequence (see \cite[Theorem~15.7]{8}),
 that for $n=2$ and 3, the natural map $\io_n:\Om_n(G)\to H_n(G)$ is a
 group isomorphism.
\end{definition}

\begin{definition}
 (\emph{the Pontryagin product} $\pp{\mu,\la}$
 in the homology group $H_2(G)$)
\smallskip

\noindent
{(a)}
Suppose that two elements $\mu,\la\in G$ commute.
Then there is a natural homomorphism $\psi:\Z\times\Z\to G$
 defined by $\psi(1,0)=\mu$, $\psi(0,1)=\la$.
The homomorphism $\psi$ induces a representation
 $\rho:\pi_1(S^1\times S^1)\to G$ and a continuous map
 $f:S^1\times S^1\to K(G,1)$.
So, we get an element $[S^1\times S^1,f]\in\Om_2(G)\cong H_2(G)$.
\smallskip

\noindent
{(b)}
The element $\pp{\mu,\la}=[S^1\times S^1,f]\in H_2(G)$
 is called \emph{the Pontryagin product}.
It is well-defined and satisfies the relation
 $\pp{\mu,\la_1}+\pp{\mu,\la_2}=\pp{\mu,\la_1\la_2}$
 when $\la_1$ and $\la_2$ both commute with $\mu$ in $G$,
 see \cite{1}.
\end{definition}

We shall use repeatedly the following identification:
Let $W$ be connected manifold with base point $x_0$.
There is a one-to-one correspondence between
 representations $\rho:\pi_1(W,x_0)\to G$
 and homotopy classes of continuous pointed maps $f:W\to K(G,1)$,
 see \cite[Chapter~6, Theorem~6.39(ii)]{8}.


\subsection{Necessity in Theorem~1.4}

\begin{lemma}
Let $G$ be a group with a meridional system
 $\vmu=(\mu_1,\ldots,\mu_r)\in G^r$.
If a system $\vla=(\la_1,\ldots,\la_r)\in G^r$ is
 weakly realizable for the pair $(G,\vmu)$, then
 condition~(i) of Theorem~1.4 holds.
\end{lemma}
\begin{proof}
Let curves $m_i,l_i$ be associated to the elements
 $\mu_i,\la_i$ by Definition~1.3b, $i=1,\ldots,r$.
Then $m_i,l_i$ lie on the boundary of a small tubular neighbourhood
 $T(L_i)$ of $L_i\subset L$.
Since $\pi_1(\bd T(L_i))\cong\Z\oplus\Z$, then the corresponding
 loops $m_i,l_i\in\pi_1(S^3-L)$ commute.
Hence their images $\rho(m_i)=\mu_i$, $\rho(l_i)=\la_i$ commute
 in $G$.
\end{proof}

\begin{lemma}
Under the conditions of Theorem~1.4, for each $i=1,\ldots,r$,
 fix an element $\la_i\in G$ commuting with $\mu_i$.
Suppose that there exists an oriented compact \emph{connected}
 3-manifold $M$ with
 a representation $\rho:\pi_1(M)\to G$ such that
$$\bd M=\sqcup_{i=1}^r (S_i^1\times S_i^1),\;
  \rho|_{\bd M}(\pt\times S_i^1)=\mu_i,\;
  \rho|_{\bd M}(S_i^1\times \pt)=\la_i,\; i=1,\dots,r.$$

\noindent
Then the sum of the Pontryagin products $\sum_{i=1}^r \pp{\mu_i,\la_i}$
 vanishes in the group $H_2(G)$.
\end{lemma}
\begin{proof}
Let $f:M\to K(G,1)$ be a continuous map corresponding to
 the homomorphism $\rho:\pi_1(M)\to G$.
We have $[\bd M,f|_{\bd M}]=
 \sum_{i=1}^r [S_i^1\times S_i^1,f|_{S_i^1\times S_i^1}]=
 \sum_{i=1}^r \pp{\mu_i,\la_i}$ in the group $H_2(G)$.
Since the disjoint union $\sqcup_{i=1}^r (S_i^1\times S_i^1)$
 is bounded by an oriented compact connected 3-manifold $M$,
 then $[\bd M,f|_{\bd M}]=0$, i.e.
 $\sum_{i=1}^r\pp{\mu_i,\la_i}=0$.
\end{proof}

\begin{lemma}
Let $G$ be a group with a meridional system
 $\vmu=(\mu_1,\ldots,\mu_r)\in G^r$.
If a system $\vla=(\la_1,\ldots,\la_r)\in G^r$ is
 weakly realizable for the pair $(G,\vmu)$, then
 condition~(ii) of Theorem~1.4 holds.
\end{lemma}
\begin{proof}
Let $L\subset S^3$ be a link weakly realizing the given system
 $\vla\in G^r$, see Definition~1.3b.
Let $T(L)\subset S^3$ be a sufficiently small tubular
 neighbourhood of $L$.
Then Lemma~3.5 follows from Lemma~3.4
 for $M=S^3-\inp T(L)$.
\end{proof}

Necessity in Theorem~1.4 follows directly from Lemmas~3.3 and~3.5.


\subsection{Sufficiency in Theorem~1.4}

Lemmas~3.3 and 3.5 motivate the following definition.

\begin{definition}
 (\emph{an algebraic triple} $(G,\vmu,\vla)$)
\smallskip

\noindent
Let $G$ be a group with a meridional system $\vmu\in G^r$.
Let $\vla=(\la_1,\dots,\la_r)\in G^r$ be a system such that
 conditions~(i) and (ii) of Theorem~1.4 hold.
Then $(G,\vmu,\vla)$ is called \emph{an algebraic triple}.
\end{definition}

Let $G$ be a group with a meridional system $\vmu\in G^r$.
By Definition~3.6 and Lemmas~3.3, 3.5
 if a system $\vla\in G^r$ is weakly realizable for the pair $(G,\vmu)$,
 then the triple $(G,\vmu,\vla)$ is algebraic.
Our purpose is to prove the converse.

\begin{definition}
 (\emph{a geometric triple} $(M,\rho,f)$)
\smallskip

\noindent
(a)
Let $M$ be a connected oriented compact 3-manifold
 with $\bd M=\sqcup_{i=1}^r (S_i^1\times S_i^1)$.
Suppose that there is a surjective homomorphism
 $\rho:\pi_1(M)\to G$ such that
\smallskip

\noindent
$\bu$
for \emph{the meridian} $m_i=\pt\times S_i^1\subset\bd M$,
 we have $\rho(m_i)=\mu_i$;
\smallskip

\noindent
$\bu$
for \emph{the longitude} $l_i=S_i^1\times\pt\subset\bd M$,
 we have $\rho(l_i)=\la_i$ for each $i=1,\dots,r$.
\medskip

\noindent
(b)
Take a continuous map $f:M\to K(G,1)$ associated to
 $\rho$.
Then $(M,\rho,f)$ is said to be \emph{a geometric triple}
 corresponding to the algebraic triple $(G,\vmu,\vla)$.
\end{definition}

\begin{lemma}
For any algebraic triple $(G,\vmu,\vla)$, there is
 a corresponding geometric triple $(M,\rho,f)$.
\end{lemma}
\begin{proof}
By Definition~3.6 the Pontryagin products $\pp{\mu_i,\la_i}\in H_2(G)$
 are well-defined.
We have representations
 $\rho|_{S_i^1\times S_i^1}:\pi_1(S_i^1\times S_i^1)\to G$.
There are continuous maps
 $f_i:S_i^1\times S_i^1\to K(G,1)$ such that
 $\io_2([S_i^1\times S_i^1,f_i])=\pp{\mu_i,\la_i}$.
\smallskip

By condition~(ii) of Theorem~1.4 the element
 $[\sqcup_{i=1}^r(S_i^1\times S_i^1),\sqcup_{i=1}^r f_i]=
 \io_2^{-1}(\sum_{i=1}^r \pp{\mu_i,\la_i})$ vanishes in $\Om_2(G)$.
By Definition~3.1 there are a 3-manifold $M$
 and a continuous map $f:M\to K(G,1)$ extending  the maps $f_i$.
We can, if necessary, add 1-handles to $M$ to make it connected
 and add connected sums of $S^2\times S^1$
 to make the homomorphism corresponding to $f$ surjective.
\end{proof}

\begin{lemma}
Let $(G,\vmu,\vla)$ be an algebraic triple and let
 $(M,\rho,f)$ be a corresponding geometric triple.
Denote by $W$ the closed 3-manifold
 $M\cup (S^1\times(\sqcup_{i=1}^r D_i^2))$, where
 $\pt\times\bd D_i^2$ is glued to $\pt\times S_i^1$ in $\bd M$.
Fix a point $q_i\in\bd D_i^2$, $i=1,\dots,r$.
\smallskip

\noindent
{(a)}
Any closed curve $\ga\subset M$ is isotopic, in $W$,
 to a curve $\ga'\subset M$ with $\rho(\ga')=e$.
\medskip

\noindent
{(b)}
There exists an integral surgery carrying $W$ to $S^3$ in such a way that
\smallskip

\noindent
$\bu$
the link $S^1\times\cup_{i=1}^r\{q_i\}\subset
 S^1\times(\sqcup_{i=1}^r D_i^2)\subset W$
 maps to a link $L'=\cup_{i=1}^r L'_i\subset S^3$;
\smallskip

\noindent
$\bu$
$\rho:\pi_1(M)\to G$ turns into
 a surjective homomorphism
 $\rho':\pi_1(S^3-L')\to G$;
\smallskip

\noindent
$\bu$
for the meridian $m'_i$ and a longitude $l'_i$ of $L'_i$, one has
 $\rho'(m'_i)=\mu_i$, $\rho'(l'_i)=\la_i\mu_i^{a_i}$
 for some integer $a_i$.
\end{lemma}
\begin{proof}
The proof of {(a)} is completely analogous to
 \cite[Claim on p.~140]{5}.
\smallskip

\noindent
{(b)}
Any oriented closed 3-manifold $W$ can be obtained
 from $S^3$ by an integral surgery \cite[Chapter~9, Section~I]{6}.
Using~{(a)} we can perform such a surgery along curves
 $\ga\subset M$ with $\rho(\ga)=e$.
Hence we get a link $L'\subset S^3$ and a surjective homomorphism
 $\rho':\pi_1(S^3-L')\to G$.
The meridian $m_i$ always maps to the meridian $m'_i$, i.e.
 $\rho'(m'_i)=\rho(m_i)=\mu_i$.
But the longitude $l_i$ maps to some longitude
 $l'_i m_i^{-a_i}$, hence $\rho'(l'_i)=\la_i\mu_i^{a_i}$, where
 $a_i\in\Z$.
\end{proof}

\begin{proposition}
Let $G$ be a group with a meridional system $\vmu\in G^r$.
A system $\vla$ is weakly realizable for $(G,\vmu)$
 if the triple $(G,\vmu,\vla)$ is algebraic.
\end{proposition}
\begin{proof}
Suppose that $(G,\vmu,\vla)$ is an algebraic triple.
Let $(M,\rho,f)$ be a corresponding geometric triple from
 Lemma~3.8.
By forming the connected sum of $M$ with copies of $S^1\times S^2$
 we add free generators to $\pi_1(M)$ and hence can arrange that
 the representation $\rho:\pi_1(M)\to G$ is surjective.
Apply Lemma~3.9b to the manifold $M$ and the homomorphism $\rho$.
\end{proof}

Sufficiency in Theorem~1.4 follows directly from Proposition~3.10.


\section{Johnson-Livingston Products}


\subsection{Definitions of Johnson-Livingston products}

Firstly, we introduce the usual Johnson-Livingston product
 in Definition~4.1.
The extended Johnson-Livingston product will appear in
 Definition~4.4.

\begin{definition}
 (\emph{the Johnson-Livingston product}
 $\jl{\mu,\la}\in Q(G)$)
\smallskip

\noindent
(a)
Suppose that the Pontryagin product $\pp{\mu,\la}$ of two commuting elements
 $\mu,\la\in G$ vanishes in $H_2(G)$.
Since $H_2(G)\cong \Om_2(G)$, there is an oriented compact
 connected 3-manifold $M$ with a surjective homomorphism
 $\rho:\pi_1(M)\to G$ such that
$$\bd M=S^1\times S^1,\quad
  \rho|_{\bd M}(\pt\times S^1)=\mu\quad \mbox{  and  }\quad
  \rho|_{\bd M}(S^1\times\pt)=\la.$$

\noindent
(b)
Suppose that $\la\in G'$.
Form the closed manifold $U=M\cup(D^2\times S^1)$,
 where $\bd M$ is identified with $\bd(D^2\times S^1)$
 in such a way that the meridian
 $\bd D^2\times\pt\subset D^2\times S^1$ is glued to
 $S^1\times\pt\subset\bd M$.
This gluing kills the longitude $S^1\times\pt$, which
 is in the kernel of $\hat\rho=\pr\circ\rho:\pi_1(M)\to G/G'$.
By the Seifert-van-Kampen theorem,
 $\hat\rho=\pr\circ\rho$ extends to a
 homomorphism $\ti\rho:\pi_1(U)\to G/G'$.
We get the corresponding map $\ti f:U\to K(G/G',1)$ and
 an element $[U,\ti f]\in\Om_3(G/G')\cong H_3(G/G')$.
\smallskip

\noindent
(c)
Let $q:H_3(G/G')\to Q(G)=H_3(G/G')/\pr_*(H_3(G))$ 
 be the projection map.
The element
 $\jl{\mu,\la}=q\circ\io_3([U,\ti f])\in Q(G)$
 is \emph{the Johnson-Livingston product}.
\end{definition}

Recall that the Pontryagin product is additive \cite{1}.
Then one gets a well-defined homomorphism
 $\te_{\mu}:Z(\mu)\to H_2(G)$, $\te_{\mu}(\la)=\pp{\mu,\la}$.
Here $Z(\mu)\subset G$ is the centralizer subgroup of $\mu\in G$.
Denote by $P(G,\mu)\subset Z(\mu)$ the kernel of $\te_{\mu}$.

\begin{proposition} \emph{\cite[Section~1]{5}}
For $\la\in P(G,\mu)$, the Johnson-Livingston product
 $\jl{\mu,\la}\in Q(G)$ is well-defined.
Moreover, $\jl{\mu,\la}+\jl{\mu,\tau}=\jl{\mu,\la\cdot\tau}$
 when defined.
Hence, the map $\chi_{\mu}: P(G,\mu)\to Q(G)$ is a group homomorphism.
\end{proposition}

We introduce the notion of a geometric pentad which will be used
 in the definition of the extended Johnson-Livingston product
 below.

\begin{definition}
 (\emph{a geometric pentad} $(U,F,M,\rho,f)$)
\noindent
\smallskip

\noindent
Let $(G,\vmu,\vla)$ be an algebraic triple,
 $(M,\rho,f)$ be a corresponding geometric triple.
Take a connected oriented surface $F$ with boundary
 $\bd F=\sqcup_{i=1}^r S_i^1$.
Let $p_i\in S_i^1\subset \bd F$ be a point, $i=1,\dots,r$.
Form the closed manifold $U=M\cup(F\times S^1)$,
 where $\bd F\times S^1$ is identified with $\bd M$ in such a way that
\smallskip

\noindent
$\bu$
\emph{the meridian} $m_i=\pt\times S_i^1\subset\bd M$ is glued to
 $\{p_i\}\times S^1\subset\bd F\times S^1$;
\smallskip

\noindent
$\bu$
\emph{the longitude} $l_i=S_i^1\times\pt\subset\bd M$ is glued to
 $S_i^1\times\pt\subset\bd F\times S^1$
\smallskip

\noindent
 for each $i=1,\dots,r$.
The pentad $(U,F,M,\rho,f)$ with all the above properties
 is called \emph{a geometric pentad} corresponding to
 the algebraic triple $(G,\vmu,\vla)$.
\end{definition}

In order to extend the Johnson-Livingston product, we shall
 suppose from now on that the elements $\mu_i$ of the meridional
 system $\vmu$ are conjugate to one another in $G$.
This implies that the abelianization $G/G'$ is a cyclic group.
We shall denote by $n$ its order.
In particular, $n=0$ when $G/G'\cong\Z$ and $n=1$ when $G=G'$.

\begin{definition}
 (\emph{the extended Johnson-Livingston product}
 $\jl{\vmu,\vla}\in Q(G)$)
\smallskip

\noindent
Let $(G,\vmu,\vla)$ be an algebraic triple.
Assume that $\vla\in (G')^r$.
Let $(U,F,M,\rho,f)$ be a corresponding geometric pentad.
We shall show in Lemma~4.13 that the homomorphism
 $\hat\rho=\pr\circ\rho:\pi_1(M)\to G/G'$ extends to
 a homomorphism $\ti\rho:\pi_1(U)\to G/G'$.
We get the corresponding map $\ti f:U\to K(G/G',1)$
 and an element $[U,\ti f]\in\Om_3(G/G')$.
Its image $\jl{\vmu,\vla}=q\circ\io_3([U,\ti f])$ in $Q(G)$
 is well-defined and is called
 \emph{the extended Johnson-Livingston product},
 see Theorem~4.17.
\end{definition}

\noindent
In \cite[Appendix 3]{5}, there is
 an example of $G$ such that $Q(G)$ is non-trivial.


\subsection{Multi-connected sums}

\begin{definition}
 (\emph{a multi-connected sum} of manifolds)
\smallskip

\noindent
{(a)}
Let $(G,\vmu,\vla)$ and $(G,\vmu,\vta)$ be algebraic triples.
Assume that condition~(iii) of Theorem~1.5 holds for $\vla$ and $\vta$.
Let $(U_M,F_M,M,\rho_M,f_M)$ and $(U_N,F_N,N,\rho_N,f_N)$ be
 geometric pentads corresponding to $(G,\vmu,\vla)$ and
 $(G,\vmu,\vta)$, respectively.
\smallskip

\noindent
{(b)}
Denote by $B_i$ and $C_i$ sufficiently small annular
 neighbourhoods of the meridians $\pt\times S_i^1\subset\bd M$ and
 $\pt\times S_i^1\subset\bd N$, respectively, $i=1,\dots,r$.
Let $A=\sqcup_{i=1}^r A_i$ be the disjoint union of $r$
 annuli $A_i$.
\smallskip

\noindent
{(c)}
Set $M\#_r N=M\cup (A\times[0,1])\cup N$, where
 $A_i\times\{0\}$ is identified with $B_i\subset\bd M$,
 $A_i\times\{1\}$ is identified with $C_i\subset\bd N$,
 $i=1,\dots,r$.
The manifold $M\#_r N$ is called \emph{a multi-connected sum}
 of the manifolds $M,N$, see Fig.~3.
\end{definition}

\begin{figure}[h]
\includegraphics{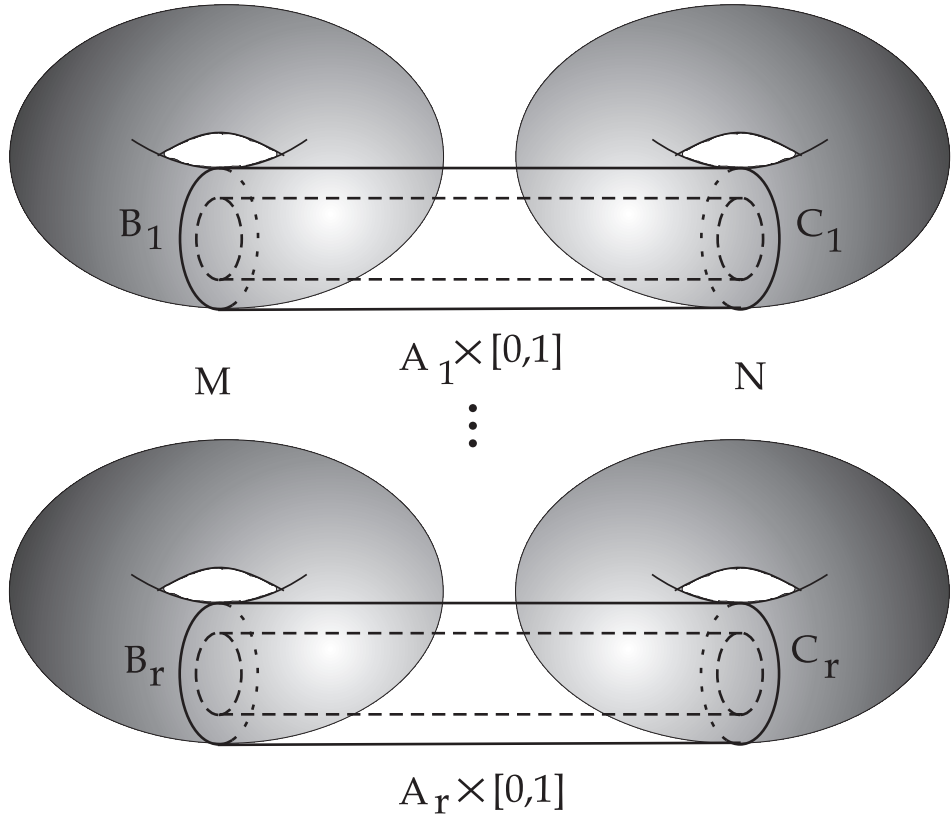}
\caption{
A multi-connected sum $M\#_r N$ of 3-dimensional manifolds.}
\end{figure}

\begin{definition}
 (\emph{a multi-band sum} of surfaces)
\smallskip

\noindent
Denote by $I$ the disjoint union $\sqcup_{i=1}^r I_i$ of $r$
 segments. Set $b_i=B_i\cap F$, $c_i=C_i\cap F$.
Set $F:=F_M\#_r F_N=F_M\cup(I\times[0,1])\cup F_N$, where
 $I_i\times\{0\}$ is identified with $b_i\subset\bd F_M$,
 $I_i\times\{1\}$ is identified with $c_i\subset\bd F_N$.
The surface $F:=F_M\#_r F_N$ is \emph{a multi-band sum}
 of the surfaces $F_M$ and $F_N$.
\end{definition}

\begin{definition}
 (\emph{a multi-connected sum} of geometric pentads)
\smallskip

\noindent
By using $M\#_r N$ and the surface $F$
 one can construct a closed manifold $W$.
By Lemma~4.8a below the homomorphisms
 $\rho_M:\pi_1(M)\to G$, $\rho_N:\pi_1(N)\to G$ and
 the continuous maps $f_M:M\to K(G,1)$, $f_N:N\to K(G,1)$
 can be extended to a homomorphism $\rho:\pi_1(M\#_rN)\to G$
 and a continuous map $f:M\#_r N\to K(G,1)$.
The pentad $(W,F,M\#_r N,\rho,f)$ is said to be
 \emph{a multi-connected sum} of pentads.
\end{definition}

\begin{lemma}
The geometric pentad $(W,F,M\#_r N,\rho,f)$ corresponds to
 the algebraic triple $(G,\vmu,\vla\cdot\vta)$.
In more details,
\smallskip

\noindent
{(a)}
the homomorphisms $\rho_M:\pi_1(M)\to G$,
 $\rho_N:\pi_1(N)\to G$ and
 the corresponding maps $f_M:M\to K(G,1)$, $f_N:N\to K(G,1)$
 can be extended to a homomorphism $\rho:\pi_1(M\#_rN)\to G$
 and a continuous map $f:M\#_r N\to K(G,1)$;
\smallskip

\noindent
{(b)}
we have $\bd(M\#_rN)=\sqcup_{i=1}^r S_i^1\times S_i^1$,
 $\rho(m_i)=\mu_i$, $\rho(l_i)=\la_i\tau_i$, where\\
 $m_i=\pt\times S_i^1\subset\bd(M\#_rN)$,
 $l_i=S_i^1\times\pt \subset\bd(M\#_rN)$, $i=1,\dots,r$.
\end{lemma}
\begin{proof}
A Seifert-Van-Kampen argument taking into account the arcs
 joining the base point to boundary components of $M$ and $N$
 shows that $\rho_M$ and $\rho_N$ can be extended to $\rho$
 with the required properties, see Fig.~3.
\end{proof}

\begin{definition}
 (\emph{a multi-connected sum} of links)
\smallskip

\noindent
A link $L=\cup_{i=1}^r L_i\subset S^3$ is called
 \emph{a multi-connected sum} of links
 $K=\cup_{i=1}^r K_i\subset S^3$ and $J=\cup_{i=1}^r J_i\subset S^3$,
 if there exist a two-sided 2-sphere $S\subset S^3$ and arcs $I_i\subset S$,
 such that $K_i\cap J_i=I_i$, $(K_i\cup J_i)-I_i=L_i$
 for each $i=1,\ldots,r$, and $K$ lies inside $S$, $J$ lies outside $S$.
In other words, we simultaneously make $r$ usual connected sums
 on different components of $K$ and $J$.
We shall denote a multi-connected sum by $L=K\#_r J$, see Fig.~4.
\end{definition}

\begin{figure}[h]
\includegraphics{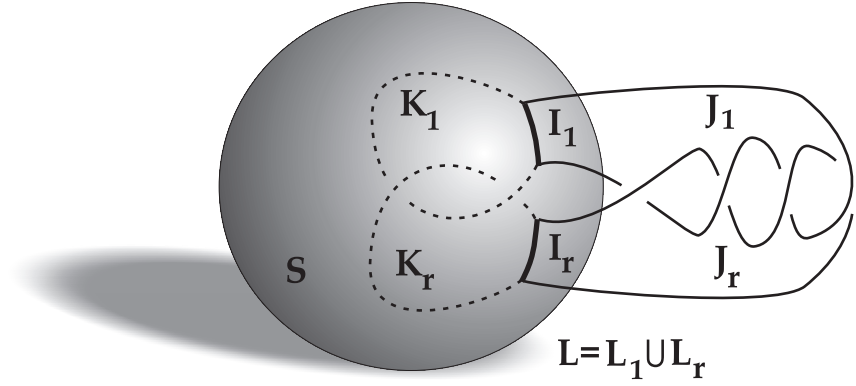}
\caption{
A multi-connected sum $L=K\#_J N$ of links.}
\end{figure}

If $r=1$, then $L=K\#_1 J$ is the usual \emph{connected sum} of
 knots.
The following lemma is geometrically obvious.

\begin{lemma}
Let $L=\cup_{i=1}^r L_i$ and $L'=\cup_{i=1}^r L'_i$ be links in $S^3$.
If $(l_1,\dots,l_r)$ and $(l'_1,\dots,l'_r)$ are systems
 of preferred longitudes for the link $L$ and $L'$, respectively,
 then the usual connected sums $(l_1\#_1 l'_1,\dots,l_r\#_1 l'_r)$
 form a preferred system of longitudes for the multi-connected
 sum $L\#_r L'$.
\end{lemma}

Actually, if $\bd F=\cup_{i=1}^r l_i$ and
 $\bd F'=\cup_{i=1}^r l'_i$, then
 $\bd (F\#_r F')=\cup_{i=1}^r(l_i\#_1 l'_i)$.

\begin{proposition}
{(a)}
The set $R(G,\vmu)$ is a subgroup of $G^r$.
The trivial element of $R(G,\vmu)$ is represented by
 an algebraically split link.
If links $K$ and $J$ realize systems $\vla,\vta\in G^r$,
 respectively, then a multi-connected sum $K\#_r J$ realizes
 $\vla\cdot\vta$.
\smallskip

\noindent
{(b)}
Let $(\mu_1,\dots,\mu_r)$ be a meridional system of a group $G$
 and suppose that $(\la_1,\dots,\la_r)$ is realized
 by a link $L\subset S^3$.
Then, for any $\mu_{r+1}\in G$,
 $(\mu_1,\dots,\mu_r,\mu_{r+1})$ is a meridional system and
 $(\la_1,\dots,\la_r,e)$ is realized by the link
 $L\sqcup O\subset S^3$, i.e.
 the link $L$ with a split off trivial component added to it.
\end{proposition}
\begin{proof}
Let $K=\cup_{i=1}^r K_i\subset S^3$ and $J=\cup_{i=1}^r J_i\subset S^3$
 be two links realizing systems $\vla,\vta\in R(G,\vmu)$.
We are going to construct a link
 $L=\cup_{i=1}^r L_i\subset S^3$ realizing the system
 $\vla\cdot\vta=(\la_1\tau_1,\dots,\la_r\tau_r)$.
Take a multi-connected sum $L=K\#_r J$.
Let $B,C\subset S^3$ be the internal and external 3-balls
 bounded by the sphere $S$ from Definition~4.9.
Then the following diagram is commutative.

$$\begin{CD}
\pi_1(S-\cup_{i=1}^r\bd I_i) @>>> \pi_1(B-\cup_{i=1}^r(K_i-\inp I_i))\\
@VVV @VVV \\
\pi_1(C-\cup_{i=1}^r(J_i-\inp I_i)) @>>> \pi_1(S^3-L)
\end{CD}$$
\smallskip

By the Seifert-van-Kampen theorem, we may extend the given
 surjective representations $\rho_K:\pi_1(S^3-K)\to G$ and
 $\rho_J:\pi_1(S^3-J)\to G$ to a surjective representation
 $\rho_L:\pi_1(S^3-L)\to G$ with peripheral information
 specified as required by Lemma~4.10.
\smallskip

To realize the inverse of a system $\vla$ in $R(G,\vmu)$,
 realize $\vla$ by a link $L$ and reverse the orientations on the
 link components.
Take a ribbon link $L$ which realizes some element
 $\vla\in R(G,\vmu)$ by using Proposition~2.3.
The multi-connected sum of $L$ and its inverse
 is an algebraically split link which realizes the trivial
 element.
\smallskip

\noindent
The item {(b)} is geometrically obvious.
\end{proof}


\subsection{Lifting continous maps and homomorphisms}

\begin{lemma}
Let $H$ be a group, $U_1,U_2,W$ be compact $n$-dimensional
 manifolds with boundary.
Take orientation preserving diffeomorphisms
 $g_1:\bd W\to\bd U_1$ and $g_2:\bd W\to\bd U_2$.
Let $f_1:U_1\to K(H,1)$, $f_2:U_2\to K(H,1)$, $f:W\to K(H,1)$
 be continuous maps such that
 the following diagrams are commutative:
$$\begin{CD}
\bd W @>{f|_{\bd W}}>> K(H,1) @. \hspace{2cm} @.
\bd W @>{f|_{\bd W}}>> K(H,1)\\
@V{g_1}VV  @|  @.  @V{g_2}VV  @| \\
\bd U_1 @>{f_1|_{\bd U_1}}>> K(H,1), @. \hspace{2cm} @.
\bd U_2 @>{f_2|_{\bd U_2}}>> K(H,1).\\
\end{CD}$$

\noindent
Then $f_1,f_2$, and $f$ extend to continuous maps
$$\ti f_1:U_1\cup_{g_1}(-W)\to K(H,1)\quad\mbox{ and }\quad
  \ti f_2:U_2\cup_{g_2}(-W)\to K(H,1)$$
 such that in the group $\Om_n(H)$ we have
$$[U_1\cup_{g_1}(-W),\ti f_1]-[U_2\cup_{g_2}(-W),\ti f_2]=
  [U_1\cup_{g_3}(-U_2),\ti f_3],$$
 where $g_3=g_1\circ g_2^{-1}: \bd U_2\to\bd W\to\bd U_1$
 and $\ti f_3$ extends $f_1$ and $f_2$.
\end{lemma}
\begin{proof}
It is clear that the maps $f_1,f_2$, and $f$ extend.
Consider the $(n+1)$-dimensional manifold with corners
$$X=(U_1\cup_{g_1}(-W))\times[0,1]\cup_{h_1}
    (W\times[0,1])\cup_{h_2}(U_2\cup_{g_2}(-W))\times[0,1],$$
 where $h_1$ identifies $W\times\{0\}$ with
 $(-W)\times\{1\}\subset(U_1\cup_{g_1}(-W))\times\{1\}$ and
 $h_2$ identifies $W\times\{1\}$ with
 $(-W)\times\{0\}\subset(U_2\cup_{g_2}(-W))\times\{0\}$.
The maps $\ti f_1$ and $\ti f_2$ extend to a map
 $\ti f:X\to K(H,1)$.
\smallskip

The boundary of $X$ has three components diffeomorphic
 to $U_1\cup_{g_1}(-W)$, $(-U_2)\cup_{g_2}W$, and
 $(-U_1)\cup(\bd W\times[0,1])\cup U_2$ respectively
 with corresponding maps to $K(H,1)$.
This proves the equality in $\Om_n(H)$
 since $[X,\ti f]=0$ in $\Om_{n+1}(H)$.
\end{proof}

\begin{lemma}
Let $(U,F,M,\rho_M,f_M)$ be a geometric pentad corresponding
 to an algebraic triple $(G,\vmu,\vla)$.
If $\vla\in (G')^r$, then
 the homomorphism $\hat\rho_M=\pr\circ\rho_M:\pi_1(M)\to G\to G/G'$
 can be extended to a homomorphism $\ti\rho_M:\pi_1(U)\to G/G'$.
\end{lemma}
\begin{proof}
A Mayer-Vietoris sequence argument shows that it suffices
 to find a homomorphism $\phi_M:H_1(F\times S^1)\to G/G'$
 such that the diagram is commutative:
$$\begin{CD}
H_1(\bd F\times S^1) @>>> H_1(F\times S^1)\\
@VVV @V{\phi_M}VV \\
H_1(M) @>{\hat\rho_M}>> G/G'
\end{CD}$$
One has $H_1(F\times S^1)\cong
 H_1(F)\otimes H_0(S^1)\oplus H_0(F)\otimes H_1(S^1)$.
Define the restriction of $\phi_M$ to $H_1(F)\otimes H_0(S^1)$
 to be the zero map.
Let the generator of $H_0(F)\otimes H_1(S^1)$ map under
 $\phi_M$ to $[\mu_1]$.
The diagram above commutes since
 $\hat\rho_M(l_i)=0$ and $\hat\rho_M(m_i)=[\mu_1]$ as all
 $\mu_i$ are conjugate to one another.
\end{proof}

\begin{lemma}
Under the conditions of Lemma~4.13
 suppose that $\la_1=\cdots=\la_r=e$ in $G$.
Then the homomorphism $\ti\rho:\pi_1(U)\to G/G'$
 from Lemma~4.13 can be lifted to
 a homomorphism $\bar\rho:\pi_1(U)\to G$.
\end{lemma}
\begin{proof}
Set $Y_1=M$ and $Y_2=T\cup(F\times S^1)$, where
 $T$ is a tree connecting the components of $\bd M$
 with a base point $q\in M$ inside $M$.
Then $Y=Y_1\cap Y_2=T\cup(\sqcup_{i=1}^rS_i^1\times S_i^1)$
 is connected and
 $\pi_1(Y)\cong\Pi_{i=1}^r(\Z[l_i]\times\Z[m_i])$.
Denote by $\ga_1$ the edge of $T$, which connect
 the base point $q$ with the first component of
 $\bd F\times S^1$.
Let $g(F)$ be the genus of the surface $F$.
Then $\pi_1(\ga_1\cup (F\times S^1))$ is generated by
 $2g(F)+r+1$ generators
 $x_1,y_1,\dots,x_{g},y_{g},l_1,\dots,l_r,z$
 and the relations $l_1\dots l_r=\Pi_{j=1}^{g(F)}[x_j,y_j]$,
 $zx_j=x_jz$, $zy_j=y_jz$, $zl_i=l_iz$, $i=1,\dots,r$,
 $j=1,\dots,g(F)$.
\smallskip

The elements $l_i$ represent the longitudes,
 the letter $z=m_1$ denotes the meridian of the first component.
To get $T\cup(F\times S^1)$ let us connect
 $q$ with the $i$-th component by an arc $\ga_i$, $i=1,\dots,r$.
This adds generators $u_1,\dots,u_r$.
Define the homomorphism $\ph:\pi_1(T\cup(F\times S^1))\to G$
 by $\ph(x_j)=\ph(y_j)=\ph(l_i)=e$,
 $\ph(z)=\mu_1$, $\ph(u_i)=\de_{i_1}$, where
 $\de_{i_1}\in G$ such that $\mu_i=\de_{i_1}\mu_1\de_{i_1}^{-1}$.
The following diagram commutes
$$\begin{CD}
\pi_1(Y) @>>> \pi_1(Y_1) \\
@VVV @VV{\rho_M}V \\
\pi_1(Y_2) @>{\ph}>> G
\end{CD}\quad \mbox{ since }\quad
\begin{CD}
l_i @>>> l_i \\
@VVV @VV{\rho_M}V \\
l_i @>{\ph}>> e
\end{CD}\qquad \mbox{ and }\quad
\begin{CD}
m_i @>>> m_i \\
@VVV @VV{\rho_M}V \\
u_izu_i^{-1} @>{\ph}>> \mu_i
\end{CD}\quad.$$
\end{proof}

\begin{remark}
It is in the proof of Lemma~4.14 that we need the hypothesis
 that the elements of the meridional system be conjugate to one
 another.
Lemma~4.14 is essential for the well-definedness
 of the extended Johnson-Livingston product.
\end{remark}


\subsection{Well-definedness of the extended Johnson-Livingston
 product}

\begin{lemma}
Let $(W,F,M\#_r N,\rho,f)$ be a multi-connected sum
 of geometric pentads $(U_M,F_M,M,\rho_M,f_M)$ and
 $(U_N,F_N,N,\rho_N,f_N)$ corresponding to algebraic triples
 $(G,\vmu,\vla)$ and $(G,\vmu,\vta)$, respectively.
Assume that $\vla\in (G')^r$ and $\vta\in(G')^r$.
\smallskip

\noindent
Using Lemma~4.13 let $\ti\rho_M:\pi_1(U_M)\to G/G'$,
 $\ti\rho_N:\pi_1(U_N)\to G/G'$, and
 $\ti\rho:\pi_1(W)\to G/G'$ be the homomorphisms extending
 $\hat\rho_M,\hat\rho_N,\hat\rho$, respectively.
\smallskip

\noindent
Let $\ti f_M:U_M\to K(G/G',1)$, $\ti f_N:U_N\to K(G/G',1)$
 and $\ti f:W\to K(G/G',1)$ be the corresponding maps.
Then $[W,\ti f]=[U_M,\ti f_M]+[U_N,\ti f_N]$
 in $\Om_3(G/G')$.
\end{lemma}
\begin{proof}
Using the notations of Definition~4.5 and applying Lemma~4.12
 for the group $G/G'$, we have
$$[W,\ti f]-[U_M,\ti f_M]=
  [(A\times[0,1])\cup N\cup (I\times[0,1]\cup F_N)\times S^1,\ti f_0],$$
 where $\ti f_0$ is the restriction of $\ti f$ to
 $(A\times[0,1])\cup N\cup
  (I\times[0,1]\cup F_N)\times S^1$.
Applying Lemma~4.12 again, we get:
$$[W,\ti f]-[U_M,\ti f_M]-[U_N,\ti f_N]=
  [(A\times[0,1])\cup(I\times[0,1])\times S^1,\ti f_1],$$
 where $\ti f_1$ is the restriction of $\ti f$ to
 $(A\times[0,1])\cup(I\times[0,1])\times S^1$.
The latter manifold is a disjoint union of $r$ copies of
 $S^2\times S^1$ and any map $S^2\times S^1\to K(G,1)$ extends
 to $B^3\times S^1$.
Therefore, $[W,\ti f]-[U_M,\ti f_M]-[U_N,\ti f_N]=0$
 in $\Om_3(G/G')$.
\end{proof}

\begin{theorem}
Assume that conditions~(i),(ii),(iii) of Theorems~1.4, 1.5 hold.
\smallskip

\noindent
{(a)}
The extended Johnson-Livingston product
 $\jl{\vmu,\vla}\in Q(G)$ is well-defined.
\smallskip

\noindent
{(b)}
We have $\jl{\vmu,\vla}+\jl{\vmu,\vta}=
 \jl{\vmu,\vla\cdot\vta}$ when defined.
\end{theorem}
\begin{proof}
{\bf (a)}
Let $(U,F,M,\rho,f)$ and $(U',F',M',\rho',f')$
 be two geometric pentads representing
 the algebraic triple $(G,\vmu,\vla)$.
Then the pentad $(-U',-F',-M',\rho',f')$ represents
 $(G,\vmu,\vla^{-1})$.
The multi-connected sum of $(U,F,M,\rho,f)$ and
 $(-U',-F',-M',\rho',f')$ represents
 $(G,\vmu,\ved)$ by Lemma~4.8 and gives
 an element in $\pr_*H_3(G)$ by Lemma~4.14.
Lemma~4.16 implies that
 $[U,\ti f]-[U',\ti f']\in\pr_*H_3(G)$, hence
 $\jl{\vmu,\vla}$ is well-defined.
\smallskip

\noindent
{\bf (b)}
The additivity is a direct consequence of Lemma~4.16.
\end{proof}

The following lemma will be used in Lemma~5.8, section 5.3.

\begin{lemma}
Suppose that conditions (i),(ii) of Theorem~1.4 hold for systems
 $\vmu=(\mu_1,\dots,\mu_r)$ and $\vla=(\la_1,\dots,\la_r)$.
Assume that $\la_r=e$ in the group $G$.
Set $\vmu'=(\mu_1,\dots,\mu_{r-1})$ and $\vla'=(\la_1,\dots,\la_{r-1})$.
Then the extended Johnson-Livingston product $\{\vmu',\vla'\}\in Q(G)$
 is well-defined and equal to $\{\vmu,\vla\}$.
\end{lemma}
\begin{proof}
It is clear that $(G,\vmu',\vla')$ is an algebraic triple.
Let $(U',F',M',\rho',f')$ be a geometric pentad realizing it.
Denote by $F$ the surface obtained by removing the interior
 of a disk $D^2$ in the interior of $F'$, hence 
 $U'=(F\times S^1)\cup (D^2\times S^1)\cup M'$.
\smallskip

Consider the disjoint union $M_0=M'\sqcup(D^2\times S^1)$.
To make it connected, perform a 1-handle surgery.
The resulting 3-dimensional manifold $M$ satisfies
 $\pi_1(M)=\pi_1(M')*\Z$, so that $\rho'$
 can be extended to $\rho:\pi_1(M)\to G$ by sending
 the free generator to $\mu_r$.
Denote by $f:M\to K(G,1)$ the corresponding map.
\smallskip

Set $U=M\cup(F\times S^1)$.
It is easy to see that $(U,F,M,\rho,f)$ is a geometric pentad
 corresponding to the triple $(G,\vmu,\vla)$.
One can then extend $\rho'$ and $\rho$ using Lemma~4.10
 to $\ti\rho':\pi_1(U')\to G/G'$ and
 $\ti\rho:\pi_1(U)\to G/G'$.
Consider the associated maps
 $\ti f':U'\to K(G/G',1)$ and $\ti f:U\to K(G/G',1)$.
\smallskip

Set $W=U\cup(F\times S^1\times[0,1])\cup(-U')$, where
 $F\times S^1\times\{0\}$ is identified with
 $F\times S^1\subset U$ and
 $F\times S^1\times\{1\}$ is identified with
 $F\times S^1\subset (-U')$.
The compact 4-manifold $W$ has the boundary $\bd W=U\cup(-U')$.
Moreover, the maps $\ti f$ and $\ti f'$ extend to $W$
 since they agree on $F\times S^1$.
Hence, the equality $[U',\ti f']=[U,\ti f]$ holds in $\Om_3(G/G')$.
This shows that $\jl{\vmu,\vla}=\jl{\vmu',\vla'}$ in $Q(G)$.
\end{proof}


\section{Realizable Systems: Proof of Theorem~1.5}


\subsection{Necessity in Theorem~1.5}

Proposition~2.3 implies that the subset $R(G,\vmu)\subset G^r$
 is non-empty.
We still assume that the elements $\mu_i$ of
 the system $\vmu$ are conjugate to one another in $G$.
\smallskip

\begin{lemma}
If $\vla\in R(G,\vmu)$, then Condition (iii) of Theorem~1.5 holds.
\end{lemma}
\begin{proof}
Let $L=\cup_{i=1}^r L_i\subset S^3$ be a link realizing
 the system $\vla\in G^r$.
In other words, one has $\rho(m_i)=\mu_i$ and $\rho(l_i)=\la_i$, where
 $m_1,\dots,m_r$ are meridians of $L_1,\dots,L_r$,
 $\ov{(l_1,\dots,l_r)}$ is a preferred system of longitudes
 for the link $L$.
By Lemma~2.2 in the group $H_1(S^3-L)$ one has
 $[l_i]=\sum_{j=1}^r\al_{ij}[m_j]$,
 where $\al_{ij}=\lk(L_i,L_j)$, $i\neq j$ and
 $\al_{ii}=-\sum_{j\neq i}\lk(L_i,L_j)$.
\smallskip

The homomorphism $\hat\rho=\pr\circ\rho:\pi_1(S^3-L)\to G\to G/G'$
 factorizes through $H_1(S^3-L)$.
Hence in the quotient $G/G'$ one gets
 $[\la_i]=\sum_{j=1}^r\al_{ij}[\mu_j]$.
By the conditions of Theorem~1.5 one has $[\mu_1]=\dots=[\mu_r]$.
Since $\sum_{j=1}^r\al_{ij}=0$, then
 $[\la_i]=0$ for each $i=1,\dots,r$.
\end{proof}

\begin{lemma}
If $\vla\in R(G,\vmu)$, then Condition (iv) of Theorem~1.5 holds.
\end{lemma}
\begin{proof}
By Lemma~5.1 and Theorem~4.17a the Johnson-Livingston product
 $\jl{\vmu,\vla}$ is well-defined.
Let $L=\cup_{i=1}^r L_i\subset S^3$ be a link with a surjective
 homomorphism $\rho:\pi_1(S^3-L)\to G$ realizing
 the given system $\vla\in G^r$.
Denote by $f:S^3-L\to K(G,1)$ a continuous map corresponding
 to $\rho$.
Let $F\subset S^3$ be an oriented surface such that $F\cap L=\bd F=L$.
Consider the sphere $S^3$ as the boundary of the ball $B^4$.
Push $F$ into $B^4$ leaving $\bd F$ in $S^3$.
Take a sufficiently small regular neighbourhood $T(F)\subset B^4$.
Then the complement $V=B^4-\inp T(F)$ is an oriented compact
 4-dimensional manifold with boundary
 $\bd V=(S^3-\cup_{i=1}^r \inp T(L_i))\cup(S^1\times F)$.
\smallskip

The pentad $(\bd V,F-\inp T(L),S^3-\inp T(L),\rho,f)$ is a geometric
 pentad corresponding to $(G,\vmu,\vla)$.
Note that $H_1(V)\cong\Z$.
Denote by $j$ the inclusion $j:S^3-\inp T(L)\to V$.
Consider the homomorphism $\ti\rho:\pi_1(V)\to G/G'$
 induced by the map which sends the generator of $H_1(V)$
 to $[\mu_1]$.
This makes the following diagram commutative
$$\begin{CD}
\pi_1(S^3-\inp T(L)) @>{\pr\circ\rho}>> G/G'\\
@V{j_*}VV @| \\
\pi_1(V) @>{\ti\rho}>> G/G'
\end{CD}$$
\smallskip

Denote by $i_*$ the homomorphism induced by the inclusion
 $\bd V\subset V$.
Let $\ti f:\bd V\to K(G/G',1)$ be the map corresponding to
 $\ti\rho\circ i_*$.
The construction above shows that $[\bd V,\ti f]$
 vanishes in the group $\Om_3(G/G')$.
This implies $\jl{\vmu,\vla}=0$.
\end{proof}


\subsection{Partial realization results}

Here we shall realize some auxilary systems needed
 for sufficiency in Theorem~1.5.
For any meridional system $\vmu\in G^r$,
 we have a well-defined homomorphism
 $\te_{\vmu}:Z(\vmu)\to H_2(G)$,
 $\te_{\vmu}(\vla)=\sum\limits_{i=1}^r\pp{\mu_i,\la_i}$.
Here $Z(\vmu)\subset G^r$ is the centralizer subgroup of $\vmu\in G^r$.
Denote by $P(G,\vmu)\subset Z(\vmu)$ the kernel of $\te_{\vmu}$.
Therefore, a system $\vla\in G^r$ is weakly realizable if and
 only if $\vla\in P(G,\vmu)$.
Recall that $n$ is the order of $G/G'$.

The following lemma is a generalization of
 \cite[Theorem~3]{5}.

\begin{lemma}
Let $\vmu\in G^r$ be a meridional system of a finitely generated
 group $G$.
If Condition~(iii) of Theorem~1.5 holds and $\vla\in P(G,\vmu)$,
 then there are $a_1,\ldots,a_r\in\Z$ such that
 $\la_{\va}:=(\la_1\mu_1^{a_1},\ldots,\la_r\mu_r^{a_r})\in R(G,\vmu)$
 and $a_i\equiv 0\pmod{n}$.
\end{lemma}
\begin{proof}
By Lemma~3.9b, there is a link $L'\subset S^3$ with
 a surjective homomorphism $\rho':\pi_1(S^3-L')\to G$
 realizing a system $\vla_{\va}:=(\la_1\mu_1^{a_1},\ldots,\la_r\mu_r^{a_r})$
 with $a_1,\dots,a_r\in\Z$.
For $i=1,\dots,r$, Condition~(iii) of Theorem~1.5 says that
 $[\la_i]=0$ in $G/G'$ and Lemma~5.1 shows that
 $[\la_i]+a_i[\mu_i]=0$ in $G/G'$.
Since $[\mu_1]$ generates $G/G'$,
 we get $a_i\equiv 0\pmod{n}$.
\end{proof}

\begin{lemma}
Let $\vmu\in G^2$ be a meridional system of a finitely generated
 group $G$.
Let $L=L_1\cup L_2\subset S^3$ be a 2-component link
 with a surjective homomorphism $\rho:\pi_1(S^3-L)\to G$
 such that $\rho(m_1)=\mu_1$, $\rho(m_2)=\mu_2$, where
 $m_1,m_2$ are the meridians of $L_1,L_2$.
\smallskip

Let $b,c$ be integers such that $b+c\equiv 0\pmod{n}$.
Then there is an oriented trivial knot $J\subset S^3-L$ such that
 $\rho(J)=e$, $\lk(J,L_1)=b$, $\lk(J,L_2)=c$.
\end{lemma}
\begin{proof}
Let us produce two oriented curves $J_1,J_2\subset S^3-L$ such that
\smallskip

$\bu$ $\rho(J_1)=\mu_1^b\mu_2^c$, $\quad\lk(J_1,L_1)=b$, $\quad\lk(J_1,L_2)=c$;

$\bu$ $\rho(J_2)=\mu_1^b\mu_2^c$, $\quad\lk(J_2,L_1)=0$, $\quad\lk(J_2,L_2)=0$.
\smallskip

To construct $J_1$ choose a closed curve in $S^3-L$ homotopic to
 $m_1^bm_2^c$, perform a homotopy to make it embedded.
Since $\mu_1^b\mu_2^c\in G'$ and $G$ is normally generated by
 $\mu_1$, one gets
 $\mu_1^b\mu_2^c=\prod_{i=1}^s\xi_{i}\mu_1^{\e_i}\xi_{i}^{-1}$,
 where the $\xi_i$ are in $G$ and the $\e_i$ are integers such that
 $\sum_{i=1}^s\e_i=0$ (\cite[Claim in the proof of Theorem~2]{5}).
\smallskip

Let $\ga_i$, $i=1,\dots,s$ be closed curves in $S^3-L$ such that
 $\rho(\ga_i)=\xi_i$.
The curve $\ga=\prod_{i=1}^s\ga_i \mu_1\ga_i^{-1}$ can be homotoped
 to a simple closed curve $J_2$ in $S^3-(L\cup J_1)$
 such that $\rho(J_2)=\mu_1^b\mu_2^c$ and
 $\lk(J_2,L_1)=\lk(J_2,L_2)=0$.
Consider the connected sum of $J_1$ and $J_2$ in $S^3-L$.
It can be unknotted by a homotopy in $S^3-L$.
The resulting curve $J$ satisfies the required properties.
\end{proof}

\noindent
The following lemma is a generalization of \cite[Theorem~2]{5}
 to 2-component links.

\begin{lemma}
Let $\vmu\in G^2$ be a meridional system of a finitely generated
 group $G$.
For all integers $b,c$ such that $b+c\equiv 0\pmod{n}$,
 one has $(\mu_1^{b^2+bc},\mu_2^{bc+c^2})\in R(G,\vmu)$.
\end{lemma}
\begin{proof}
By Proposition~4.11a there is a 2-component algebraically split link
 $L=L_1\cup L_2$ realizing the system $(e,e)\in G^2$.
Let $m_1,m_2$ be the meridians of $L_1,L_2$ and denote by
 $\bar l_1,\bar l_2$ the preferred longitudes of $L_1,L_2$,
 respectively.
Since $L$ is algebraically split,
 $(\bar l_1,\bar l_2)$ is the preferred system of longitudes for $L$.
Take a homomorphism $\rho:\pi_1(S^3-L)\to G$
 such that $\rho(m_1)=\mu_1$, $\rho(m_2)=\mu_2$,
 $\rho(\bar l_1)=\rho(\bar l_2)=e$.
By Lemma~5.4 there exists an oriented trivial knot $J\subset S^3-L$
 such that $\rho(J)=e$, $\lk(J,L_1)=b$, $\lk(J,L_2)=c$.
\smallskip

The integral surgery on $J$ with framing $+1$ carries $S^3$ to
 itself and the link $L$ to a link $L'=L'_1\cup L'_2$ with a surjective
 homomorphism $\rho':\pi_1(S^3-L')\to G$.
The meridians of $L'$ satisfy $m'_1=m_1$, $m'_2=m_2$.
Let us denote by $\ov{(l'_1,l'_2)}$ the preferred system of longitudes
 of $L'$.
Our aim is to prove that $\rho'(l'_1)=\mu_1^{b^2+bc}$ and
 $\rho'(l'_1)=\mu_2^{bc+c^2}$.
\smallskip

Let $m_J$ be the meridian of $J$ and let $l_J$ be
 a longitude such that $\lk(l_J,J)=+1$.
The homology class of $l_J$ in the group
 $H_1(S^3-L-J)\cong\Z[m_1]\oplus\Z[m_2]\oplus\Z[m_J]$ is
 given in this basis by the vector
 $[l_J]=\left(\begin{array}{c}
 b \\
 c \\
 1
 \end{array}\right)$.
Similarly, $H_1(S^3-L')\cong\Z[\mu'_1]\oplus\Z[\mu'_2]$.
In these bases the attaching map of the surgery disc
 is given by
 $\left(\begin{array}{c}
 b \\
 c \\
 1
 \end{array}\right)$.\\
The inclusion $\psi:S^3-(L\cup J)\to S^3-L'$ is described in
 homology by
 $\left(\begin{array}{ccc}
 1 & 0 & -b \\
 0 & 1 & -c
 \end{array}\right)$.
\smallskip

Denote by $\bar l'_1,\bar l'_2$ the preferred longitudes of
 the components $L'_1,L'_2$, respectively.
One gets
 $[\bar l_1]=\left(\begin{array}{c}
 0 \\
 0 \\
 b
 \end{array}\right),\quad
 [\bar l_2]=\left(\begin{array}{c}
 0 \\
 0 \\
 c
 \end{array}\right)\in H_1(S^3-L-J)$;
 $\psi_*([\bar l_1])=\left(\begin{array}{c}
 -b^2 \\
 -bc
 \end{array}\right),\;
 \psi_*([\bar l_2])=\left(\begin{array}{c}
 -bc \\
 -c^2
 \end{array}\right),\quad
 [\bar l'_1]=\left(\begin{array}{c}
 0 \\
 *
 \end{array}\right),\;
 [\bar l'_2]=\left(\begin{array}{c}
 * \\
 0
 \end{array}\right)\in H_1(S^3-L')$.
\medskip

In particular, for the link $L'$,
 one gets $\lk(L'_1,L'_2)=-bc$ and
 $l'_1=m_1^{bc}\bar l'_1$, $l'_2=m_2^{bc}\bar l'_2$.
Since $\psi(\bar l_1)$ and $\psi(\bar l_2)$ are longitudes of
 the components $L'_1$ and $L'_2$, then
 $\bar l'_1=m_1^{d_1}\psi(\bar l_1)$ and
 $\bar l'_2=m_2^{d_2}\psi(\bar l_2)$ for $d_1,d_2\in\Z$.
By substituting in the vectors above one obtains $d_1=b^2$, $d_2=c^2$.
One computes
 $l'_1=m_1^{b^2+bc}\psi(l_1)$ and $l'_2=m_2^{bc+c^2}\psi(l_2)$.
Since $\rho(l_1)=\rho(l_2)=e$, then
 $\rho'(l'_1)=\mu_1^{b^2+bc}$ and $\rho'(l'_1)=\mu_2^{bc+c^2}$
 as desired.
\end{proof}

\begin{lemma}
Let $\vmu\in G^2$ be a meridional system of a finitely generated group $G$.
\smallskip

\noindent
{(a)} $(1,\mu_2^{n^2})\in R(G,(\mu_1,\mu_2))$, where
 $n$ is the order of $G/G'$.
\smallskip

\noindent
{(b)} For any integer $h\in\Z$,
 $(\mu_1^{hn},\mu_2^{-hn})\in R(G,\vmu)$.
\end{lemma}
\begin{proof}
{(a)}
Theorem~2 of \cite{5} states that in the particular case
 $r=1$ of a knot $\mu$ is realizable, i.e. there is a knot
 $K$ with a surjective homomorphism $\rho:\pi_1(S^3-K)\to G$
 such that $\rho(m)=\mu$, $\rho(l)=\mu^{n^2}$, where
 $m$ and $l$ are the meridian and the preferred longitude of $K$.
By Proposition~4.11b the link $O\sqcup K$ realizes
 $(1,\mu_2^{n^2})$.
\smallskip

\noindent
{(b)}
Lemma~5.5 with $b=h$, $c=n-h$ implies
 $(\mu_1^{hn},\mu_2^{n^2-hn})\in R(G,\vmu)$.
By Proposition~4.11a the set $R(G,\vmu)$ is a group.
Then we get\\
 $(\mu_1^{hn},\mu_2^{-hn})=
  (\mu_1^{hn},\mu_2^{n^2-hn})\cdot(1,\mu_2^{n^2})^{-1}\in R(G,\vmu)$.
\end{proof}


\subsection{Sufficiency in Theorem~1.5}

\noindent
Here we finish the proof of our main Theorem~1.5 formulated in
 subsection~1.3.
\smallskip

\begin{proposition}
Under the conditions of Lemma~5.3 there exists $h\in\Z$ such that
 $(\la_1\mu_1^{hn},\la_2,\dots,\la_r)\in R(G,\vmu)$.
\end{proposition}
\begin{proof}
By Lemma~5.3 take integers $a_1,\dots,a_r\in\Z$ such that
 $\la_{\va}=(\la_1\mu_1^{a_1},\dots,\la_r\mu_r^{a_r})\in R(G,\vmu)$ and
 $a_i\equiv 0\pmod{n}$ for each $i=1,\dots,r$.
If $n=0$, then the result is obvious.
Suppose $n\geq 1$ and write $a_i=h_i n$, $h_i\in\Z$,
 $i=1,\dots,r$.
\smallskip

By Lemma~5.6b the system
 $(\mu_1^{h_2n},\mu_1^{-h_2n})$ belongs to $R(G,(\mu_1,\mu_2))$, hence
 $(\mu_1^{h_2n},\mu_2^{-h_2n},1,\dots,1)\in R(G,\vmu)$
 by Proposition~4.11b.
Then we obtain
$$\vla'_{\va}:=
  (\la_1\mu_1^{(h_1+h_2)n},\la_2,\la_3\mu_3^{h_3n},\dots,\la_r\mu_r^{h_rn})
  =\vla_{\va}\cdot(\mu_1^{h_2n},\mu_2^{-h_2n},1,\dots,1)\in R(G,\vmu).$$
Apply the same trick to the components $\la_1\mu_1^{(h_1+h_2)n}$
 and $\la_3\mu_3^{h_3n}$ of $\vla'_{\va}$ to kill $h_3$ and so on.
Finally, we get $h\in\Z$ such that
 $(\la_1\mu_1^{hn},\la_2,\dots,\la_r)\in R(G,\vmu)$ as required.
\end{proof}

The following lemma is a simple generalization of
 \cite[Theorem~4]{5}.

\begin{lemma}
Let $\vmu\in G^r$ be a meridional system of a finitely generated
 group $G$.
Let $h$ be an integer such that $\jl{\vmu,\vmu_h}=0$,
 where $\vmu_h=(\mu_1^{hn},1,\dots,1)$.
Then $\vmu_h\in R(G,\vmu)$.
\end{lemma}
\begin{proof}
By Lemma~4.18 one has $\jl{\vmu,\vmu_h}=\jl{\mu_1,\mu_1^{hn}}\in Q(G)$.
Theorem~4 of \cite{5} says that $\mu_1^{hn}\in R(G,\mu_1)$.
Then by Proposition~4.11b one gets $\vmu_h\in R(G,\vmu)$.
\end{proof}
\medskip

\noindent
{\bf Proof of Theorem~1.5.}
Theorem~1.5 says that, for a system $\vla\in P(G,\vmu)$,
 Conditions~(iii) and (iv) are equivalent to
 realizability $\vla\in R(G,\vmu)$.
Necessity of these conditions follows from necessity in
 Theorem~1.4 and Lemmas~5.1, 5.2.
\smallskip

To prove sufficiency take by Proposition~5.7 an integer $h\in\Z$
 such that $\vla_h\in R(G,\vmu)$, where
 $\vla_h=(\la_1\mu_1^{hn},\la_2,\dots,\la_r)$.
By Lemma~5.2 (necessity of condition~(iv)) one gets
 $\jl{\vmu,\vla_h}=0$ in $Q(G)$.
By Theorem~4.17b (additivity of the extended Johnson-Livingston product)
 one has $0=\jl{\vmu,\vla_h}=\jl{\vmu,\vla}+\jl{\vmu,\vmu_h}$, where
 $\vmu_h=(\mu_1^{hn},1,\dots,1)$.

Condition~(iv) of Theorem~1.5 means that $\jl{\vmu,\vla}=0$,
 hence one obtains $\jl{\vmu,\vmu_h}=0$.
Then by Lemma~5.8 we get $\vmu_h\in R(G,\mu)$.
Since by Proposition~4.11a the set $R(G,\vmu)$ is a group,
 $\vla=\vla_h\cdot\vmu_h^{-1}\in R(G,\vmu)$ as required.
\hfill $\square$


\section{Applications}

We give below examples of groups where conditions
 (i),(ii),(iii),(iv) of Theorems~1.4 and 1.5 apply.

\begin{example}
Let $G$ be the group of a classical knot $K$ and let
 $\mu_1,\dots,\mu_r$ be meridians of $K$.
Then $\vmu=(\mu,\ldots,\mu_r)\in G^r$ is
 a meridional system of $G$ and
 the only nontrivial conditions of Theorems~1.4 and 1.5
 are conditions (i) and (iii)
 since it is well-known that $H_2(G)=0$ and
 $H_3(G/G')=H_3(\Z)=0$.
\smallskip

For $i=1,\dots,r$, let $\la_i$ denote
 the preferred longitude of $K$ that commutes with $\mu_i$.
Then all systems of the form
 $\vla=(\la^{a_1},\ldots,\la^{a_r})$,
 $a_i\in\Z$, $i=1,\dots,r$ are realizable.
\smallskip

If $K$ is a hyperbolic knot, then one shows that
 these are only realizable systems since any element
 of $G$ which commutes with $\mu_i$ is parabolic
 and belongs to the peripheral subgroup of $G$
 containing $\mu_i$.
If $K$ is a composite knot, then there are other
 realizable systems because the preferred longitude
 of a summand of $K$ commutes with the meridian of $K$
 but is not a power of the preferred longitude of $K$.
\end{example}

\begin{example}
Consider virtual knot groups as described by
 See-Goo Kim in \cite{7}.
Let $G$ be the group of a virtual knot such that
 $H_2(G)$ is cyclic of order 2.
Such a knot exists \cite[section 6.3]{7}.
Its group has the presentation:
$$\langle{ a,b \; | \; b=a^{-1} b^2 a b^{-2} a,\;
 b=[b a^{-1},a^{-1}b]^{-1} b [b a^{-1},a^{-1}b]\rangle}.$$

Let $r>0$ be an integer and let $\mu$ and $\la$
 denote the meridian and preferred longitude of the virtual knot.
Set $\vmu=(\mu,\dots,\mu)$ and $\vla=(\la,\dots,\la)$ in $G^r$.
Then $\vmu$ is a meridional system of $G$,
 conditions~(i),(iii) and (iv) are satisfied,
 in particular, $G/G'\cong\Z$ so that $H_3(G/G')=0$.
The Pontryagin product $\pp{\mu,\la}$ generates $H_2(G)$
 by \cite[Theorem 15]{7}, hence condition (ii)
 is equivalent to the condition that $r$ is even.
\end{example}

\begin{example}
Let $G$ be the group described in Appendix 3 of \cite{5}.
It has the presentation
$$\langle a,b,x,y \; | \;
 b^{-1} a^n b=[x,y],\;
 b^{-2} a^n b^2=y^{-2} x y^2,\;
 a^{-1} b a =y^{-3} x y^3,\;
 a^{-2} b a^2 =x^{-2} y x^2 \rangle,$$
 where $n$ is a positive integer.
The group $G$ is normally generated by $a$ and
 $G/G'$ is cyclic of order $n$ generated by the class of $a$.
According to \cite{5} we have $H_2(G)=H_3(G)=0$ and since
 $G/G'\cong\Z_n$, we get $Q(G)=H_3(G/G')\cong\Z_n$.
\smallskip

Let $r$ be a positive integer and take
 $\vmu=(a,\dots,a)\in G^r$ and $\vla=(a^n,\dots,a^n)\in (G')^r$.
Conditions (i), (ii) and (iii) are satisfied.
We shall compute explicitly the extended Johnson-Livingston
 product $\jl{\vmu,\vla}$ and show that it is $r$ times
 the generator of $Q(G)\cong\Z_n$.
So, condition (iv) is equivalent to $n$ divides $r$.
\smallskip

Let $L^0\subset S^3$ be the trivial link with $r$ components.
Denote by $m_i^0$ and $l_i^0$, $i=1,\dots,r$
 the meridians and preferred longitudes of $L^0$.
Let $\rho_0:\pi_1(S^3-L^0)\to G$ be the homomorphism
 sending $m_i^0$ to $a$.
Let $M$ be the connected sum of $S^3-\inp T(L^0)$ with
 three copies of $S^2\times S^1$.
We can extend $\rho_0$ to a surjective homomorphism
 $\rho:\pi_1(M)\to G$.
Let $f:M\to K(G,1)$ be a corresponding continuous map.
\smallskip

Parametrize the components of $\bd T(L^0)$ by maps
 $g_i:S^1\times S^1\to \bd T(L^0)$ in such a way that
 $g_i(\pt\times S^1)=m_i^0$ and
 $g_i(S^1\times\pt)=(m_i^0)^n l_i^0$, $i=1,\dots,r$.
The geometric triple $(M,\rho,f)$ realizes
 the algebraic triple $(G,\vmu,\vla)$.
Let $F$ denote the 2-dimensional sphere punctured with $r$ holes.
Glue $F\times S^1$ to $M$ to get
 a closed manifold $U$ as described in definition~4.3
 and extend $\rho$ to $\rho_M:\pi_1(U)\to G/G'$
 with associated map $f:U\to K(G/G',1)$.
\smallskip

The way that $\rho_M$ is cinstructed in Lemma~4.14 shows
 that we may assume that the restriction of $\ti f$ to
 $F\times S^1$ is of the form
 $\ti f(x,t)=g(t)$ for $x\in F$ and $t\in S^1$, where
 $g:S^1\to K(G/G',1)$ is a map such that
 in homotopy the generator of $\pi_1(S^1)$ is sent
 to the class of $a$ in $G/G'$.
The associated pentad computes the extended
 Johnson-Livingston product $\jl{\vmu,\vla}$.
\smallskip

Set $U_1=F\times S^1$ and $U_2=(S^2-\inp F)\times S^1$
 which is a disjoint union of $r$ solid tori.
Lemma~4.12 shows that, in $\Om_3(G/G')$,
 the element $[U,\ti f]$ is equal to $[U',\ti f']$, where
 $U'$ is a connected sum of lens spaces of type $L(n,1)$
 and $\ti f'$ is the extension of $f$ to $U'$.
This uses the fact that any element of
 the form $[S^2\times S^1,h]$ vanishes in $\Om_3(G/G')$.
Therefore,
 $[U,\ti f]=r[L(n,1),\psi]$ in $\Om_3(G/G')$, where
 $\psi:L(n,1)\to K(G/G',1)$ induces a surjective homomorphism
 $\pi_1(L(n,1))\cong\Z_n\to G/G'$.
\smallskip

The result follows since the element $[L(n,1),\psi]$
 corresponds to a generator of $H_3(G)$,
 see \cite[section 3, p.~138]{5}.
\end{example}


\section*{Acknowledgements}
The authors thank Hugh Morton, Luisa Paoluzzi and Bernard Perron
 for useful discussions.
The first author was supported by Marie Curie Fellowship 007477.


\end{document}